%% file: cnincidenceschubert.tex
\documentclass{article}
\usepackage[latin1]{inputenc}
\usepackage{amsmath,amssymb,amsthm}
\usepackage[dvips]{graphicx,epsfig}

\renewcommand{\AA}{\mathbb{A}}

\newcommand{\HH}{\mathbb{H}}

\newcommand{\NN}[0]{\mathbb{N}}

\newcommand{\CC}{\raisebox{.6ex}{${\scriptscriptstyle /}$}\hspace{-.43em}C}
\newcommand{\ZZ}{\mbox{\rm \lower0.3pt\hbox{$\angle\!\!\!$}Z}}

\newtheorem{nt}{Notation}
\newtheorem{conjecture}[nt]{Conjecture}
\newtheorem{prop}[nt]{Proposition}

\newcommand{\fd}{\ensuremath{\rightarrow}}

\newcommand{\findem}{\hfill\rule{2mm}{2mm}}

\renewcommand{\phi}{\ensuremath{\varphi}}

\newcommand{\inc}{\ensuremath{\subset}}
\newcommand{\nl}{\ \\[2mm]}

\renewcommand{\phi}{\ensuremath{\varphi}}

\newcommand{\s}{Spec\;}
\newcommand{\x}{\ensuremath{\times}}

\newcommand{\barr}{\begin{array}}
\newcommand{\earr}{\end{array}}
\newcommand{\bit}{\begin{itemize}}
\newcommand{\eit}{\end{itemize}}
\newcommand{\beq}{\begin{eqnarray*}}
\newcommand{\eeq}{\end{eqnarray*}}
\newcommand{\beqn}{\begin{eqnarray}}
\newcommand{\eeqn}{\end{eqnarray}}
\newcommand{\bconj}{\begin{conjecture}}
\newcommand{\econj}{\end{conjecture}}
\newcommand{\bcor}{\begin{coro}}
\newcommand{\ecor}{\end{coro} \noindent}
\newcommand{\ben}{\begin{enumerate}}
\newcommand{\een}{\end{enumerate} \noindent}
\newcommand{\bnot}{\begin{nt} }
\newcommand{\enot}{\end{nt} \noindent}
\newcommand{\bdefi}{\begin{defi}}
\newcommand{\edefi}{\end{defi} \noindent}
\newcommand{\bprop}{\begin{prop}} 
\newcommand{\brap}{\begin{rappel}}
\newcommand{\erap}{\end{rappel} \noindent }
\newcommand{\brq}{\begin{rem}}
\newcommand{\erq}{\end{rem} \noindent }
\newcommand{\bthm}{\begin{thm}}
\newcommand{\blm}{\begin{lm}}
\newcommand{\bex}{\begin{ex}}
\newcommand{\eex}{\end{ex}\noindent }
\newcommand{\bexo}{\begin{exo} \normalfont}
\newcommand{\eexo}{\end{exo}\noindent }

\font \tengothic=eufm10 scaled\magstep 1
\font\sevengothic=eufm7 scaled\magstep 1
\newfam\gothicfam
      \textfont\gothicfam=\tengothic
      \scriptfont\gothicfam=\sevengothic
\def\goth#1{{\fam\gothicfam #1}}

\newtheorem{coro}[nt]{Corollary} 
\newtheorem{defi}[nt]{Definition} 
\newtheorem{defprop}[nt]{Definition-Proposition} 
 
\newtheorem{ex}[nt]{Example} 
\newtheorem{exo}[nt]{Exercise} 
\newtheorem{lm}[nt]{Lemma} 
\newtheorem{rappel}[nt]{Recall} 
\newtheorem{rem}[nt]{Remark} 
\newtheorem{thm}[nt]{Theorem}

\newcommand{\elm}{\end{lm} \noindent \textit{Proof: }}  
\newcommand{\eprop}{\end{prop} \noindent \textit{Proof: }} 
\newcommand{\ethm}{\end{thm}\noindent \textit{Proof: }}

\begin{document}
\title{Incidence relations among the Schubert cells of 
equivariant Hilbert
schemes}
\date{}
\author{Laurent Evain}
\maketitle
\noindent
{\bf Abstract:} Let $\HH_{ab}(H)$ be the equivariant 
Hilbert scheme  parametrizing the 
zero dimensional
subschemes of the affine plane $k^2$, fixed under the 
one dimensional torus $T_{ab}=\{(t^{-b},t^{a}),\ t\in k^*\}$ 
and whose Hilbert
function is $H$. This Hilbert scheme 
admits a natural stratification in
Schubert cells which
 extends the notion of  Schubert cells on Grassmannians. 
However, the incidence 
relations between  the cells 
become more complicated than in the case of
Grassmannians.  In this paper, we give a necessary condition 
for the closure of a  cell to meet another cell. 
In the particular case
of Grassmannians, it coincides with the well known necessary
and sufficient  incidence condition. There is no known 
example showing that the condition wouldn't be sufficient. 

\section{Introduction}
\noindent
{\bf The problem (general wording).}\\
Fixing an algebraically closed field $k$,
the Hilbert scheme $\HH$
parametrizing the zero dimensional 
subschemes of the affine plane $\s k[x,y]$ admits a natural 
action of the two dimensional torus $k^*\x k^*$,
induced by the linear action 
$(t_1,t_2).x^{\alpha}y^{\beta}=(t_1.x)^{\alpha}(t_2.y)^{\beta}$ of $k^*\x k^*$
on $k[x,y]$. If 
$T_{ab}=\{(t^{-b},t^a),\ t\in k^*\}$ is a one dimensional sub-torus,
the closed subscheme 
$\HH_{ab}$ which  parametrizes by definition the zero dimensional 
subschemes invariant under  
the action of $T_{ab}$ is a disjoint union of irreducible subschemes 
$\HH_{ab}(H)$. Each $\HH_{ab}(H)$ can be embedded in a product of 
Grassmannians and inherits  the  stratification whose cells are 
the inverse images of the products of the Schubert cells. 
The problem is to understand the geometry of this stratification,
and in particular to describe the incidences between the cells.
\nl
{\bf Motivations}\\
Before explaining the results and the techniques, let's explain the 
motivations. Roughly speaking, Grassmannians are 
easier to study than Hilbert schemes because they 
are stratified by Schubert cells. Those stratifications 
enable for instance to compute
the intersection rings of Grassmannians, the Hilbert function
of their natural embeddings in a projective space ([GH],[M]). 
On the other hand, Hilbert schemes admit  a stratification
furnished by the theory of Grobner bases (or standard basis) 
which is fussy to describe ([Gr]). 
The equivariant Hilbert scheme $\HH_{ab}(H)$ is an object which both
looks like a Hilbert scheme and a Grassmannian, 
a sort of interface between them. It has two stratifications,
the stratification in Schubert cells introduced above and 
the Gr\"obner stratification obtained by restricting the stratification 
of $\HH$ to $\HH_{ab}(H)$. Moreover these two stratifications coincide. 
So the philosophy is to try to describe the geometry 
of $\HH_{ab}=\cup \HH_{ab}(H)$,
and then to lift the information 
on $\HH$. This is our motivation to study $\HH_{ab}(H)$.
To see some precise examples 
where the geometry of $\HH_{ab}$ can be used to describe the geometry
of $\HH$, see $[I]$, where  structures of bundles are highlighted, 
[ES] where the Betti numbers 
of the connected components of $\HH$ are computed via the study 
of the equivariant inclusion 
$\HH_{ab} \inc \HH $ for $(a,b)$ general enough, 
or [B] which describes links between the equivariant 
cohomology of $\HH_{ab}$ and the equivariant cohomology of $\HH$. 
Beyond this motivation, note 
that $\HH_{ab}$ is a classical object of study ([G\"o],[I],[Y]
for instance).
\nl
{\bf Known results and precise formulation of the problem.}\\
For some special Hilbert functions $H$, the equivariant Hilbert scheme 
$\HH_{ab}(H)$ is a Grassmannian and the cells of $\HH_{ab}(H)$
are the Schubert cells (for this reason, we will say  a
``Schubert cell'' to designate a cell on $\HH_{ab}(H)$ ).
In this simple case, the closure 
of a cell is a union of cells and there is an explicit numerical condition
to check incidence. 
This is no longer so simple when $\HH_{ab}(H)$ is
not a Grassmannian: 
Yameogo has given in [Y] an example where the closure of a cell
is not a union of cells.  So the problem splits into two different pieces:
the ``weak incidence'' problem which consists in  deciding whether the closure
of a cell meets another cell, 
and the ``strong incidence'' problem which consists in 
deciding whether the closure of a cell contains another cell. 
In [Y], Yameogo faced the weak incidence problem and gave
a  necessary 
but not sufficient condition to have a relation 
$\overline{C} \cap C' \neq \emptyset$.
In this paper, we address  the weak incidence problem too.

\nl
{\bf The results.} \\
Our main answer will 
be a necessary condition (theorem \ref{thm}) 
for weak incidence. No counter-example to sufficiency is known to 
the author.
The other propositions 
try to analyse the pertinence of the condition,
by comparing it to the previous 
results (prop. \ref{idemG}, \ref{amelJ}) 
and by testing it on small Hilbert functions (section \ref{sex}).
\\
We now formulate our results with more details. We associate with 
any Schubert cell $C$ a combinatorial data and we express the relation
$\overline C \cap C' \neq \emptyset$ 
by a combinatorial property linking $C$ and $C'$, as follows.

  \begin{figure}[h] 
     \begin{center}
        \input{escetsondual.pstex_t}
     \end{center} 
  \end{figure}
Each Schubert cell is 
determined by a staircase $E$, i.e. a subset of $\NN^2$ whose complementary
is stable under the addition of $\NN^2$: the ideal $I^E$ generated by the 
monomials whose exponent is not in $E$ is the unique monomial 
ideal of the cell. If $n$ is the cardinal of $E$, $E$ is included in 
a box of size $n\x n$ and we denote by  $E^{\nu}$ the dual staircase 
defined by the complementary of $E$ in the box (see figure). 
A system of arrows $S$ on $E$ 
is the data for each element $p$ of $E$ of an arrow $(p,f_p)$
with origin $p$, the set of arrows being compelled to compatibility 
conditions. 

  \begin{figure}[h] 
     \begin{center}
        \input{exempleincidence.pstex_t}
     \end{center} 
  \end{figure} Replacing each 
element $p$ of $E$ by $f_p$ gives a new subset of $\NN^2$ which turns out 
to be a staircase thanks to the compatibility conditions. 
Thinking of $S$ as an operator on $E$, we denote by $S(E)$ 
the staircase obtained from $E$ and $S$ by the above procedure.
The theorem says that two cells $C(E)$ and $C(F)$ corresponding 
respectively to staircases $E$ and $F$
verify $\overline {C(E)} \cap C(F) \neq \emptyset$ only if 
\ben
\item
there exists a system of arrows $S$ on $E$ such that $S(E)=F$
\item
there exists a system of arrows $S^\nu$ on $E^\nu$ such that
$S^\nu(E^{\nu})=F^{\nu}$
\een
We then compare this condition to the known conditions
and we test it on an example. In  the particular
case of Grassmannians, 1) and 2) are equivalent
and  reduce  to the well 
known necessary and sufficient incidence condition (prop. \ref{idemG}).
In general, they are not equivalent (ex. \ref{condnoneq})
but any of them imply the condition of [Y] (proposition \ref{amelJ}). 
Finally, we study  an example (section \ref{sex}) which illustrates that 
when the Hilbert function $H$ is small, we can solve the weak incidence 
problem because  our
condition
can be shown to be necessary and sufficient.
\nl
\noindent
{\bf The methods.}\\
The first step consists in putting the problem in the context of Gr\"obner
basis theory and to describe a cell $C(E)$ as the locus formed by the ideals
$I$ having  $I^E$ as initial ideal. Then we interpret the incidence relation 
$\overline{C(E)} \cap C(F) \neq \emptyset$ as the existence of 
a rational morphism from $\AA^1=\s k[t]$ to ${ C(E)}$ which can be 
extended in $\infty$ by  putting $f(\infty)=X_F \in 
C(F)$. The universal ideal
on $\HH_{ab}(H)$ restricted on  the  affine curve defines an ideal 
$I(t)\inc k[x,y,t]$. The key consists in using ideas coming from
Grobner basis  theory to  exhibit within $I(t)$ 
a set of elements $P_1,P_2,\dots$ from which we read the system of arrows 
on $E^\nu$. 
Then, using an argument of duality corresponding to the geometric notion 
of linkage of two zero dimensional schemes in a complete intersection, 
we get the system of arrows on $E$.
\nl 
I thank Michel Brion who told me about the problems that motivated 
this work and who answered several questions.

\section{The stratification on $\HH_{ab}(H)$.} \label{s1}
\subsection{Notations} \label{ssnot}
First, we keep the notations of the introduction:
$k$ is an algebraically closed field,
$a$ and $b$ are two 
relatively prime integers, 
and $T_{a,b}=\{(t^{-b},t^a), t \in k^*\}$ is a one 
dimensional subtorus of $k^*\x k^*=:T$. The torus $T$ acts on 
the Hilbert scheme $\HH$ parametrizing the zero-dimensional subschemes 
of $\s k[x,y]$, and $\HH_{ab}$ is the subscheme of $\HH$ parametrizing
the subschemes invariant under the action of
$T_{ab}\inc T$.
Alternatively, the subschemes of  $\HH_{ab}$ can be characterized 
by their ideals using degrees.  
 The degrees $d,d_x,d_y$ are defined on monomials by 
$d(x^\alpha y^\beta)=-b \alpha+ a \beta,
d_x(x^\alpha y^\beta)=\alpha ,d_y(x^\alpha y^\beta)=\beta
$. If $I$ is an ideal of $k[x,y]$, we let
$I_n:=I \cap k[x,y]_n$, where $k[x,y]_n$ denotes the vector space
generated by the monomials $m$ of degree $d(m)=n$. 
A subscheme $Z$  is in $\HH_{ab}$, 
if and only if its ideal  is quasi-homogeneous with respect to $d$, 
ie. $I(Z)=\oplus_{n\geq 0} I(Z)_n$.
\\
We order the monomials 
of $k[x,y]$ by the rule $m_1<m_2$ if $d(m_1)<d(m_2)$ or ($d(m_1)=d(m_2)$ 
and $d_y(m_1)<d_y(m_2)$).
\\
The scheme $\HH_{ab}$ is not connected but the connected components 
are determined by a Hilbert function. 
By semi-continuity, if a subscheme $Z' \in \HH_{ab}$ is a 
specialization of $Z\in  \HH_{ab}$, then the codimensions 
of $I_n(Z)$ and $I_n(Z')$ in $k[x,y]_n$ verify 
$codim I_n(Z) \geq codim I_n(Z')$. 
But $Z$ and $Z'$ 
have the same length
$l=\sum _{n\geq 0}codim I_n(Z)= \sum _{n\geq 0}codim I_n(Z')$.
It follows that the 
sequence $H(Z)=(h_0,h_1,h_2 \dots )$ 
where $h_i=codim I_n(Z)$ is constant 
on the connected components 
of $\HH_{ab}$ and that $h_n=0$ for $n$ big enough. 
If $H=(h_0, \dots,h_r,0,0,\dots )$ is any sequence, 
we denote by $\HH_{ab}(H)$ 
the closed subscheme of $\HH$ parametrizing the schemes 
$Z$ verifying $H(Z)=H$. One can verify that $\HH_{ab}(H)$ is  an 
irreducible connected  component of $\HH_{ab}$
(though we won't use it in the sequel).
\\
Recall that a staircase is a subset of $\NN^2$ whose complementary
is stable by addition of $\NN^2$. Staircases will be used 
to parametrize the stratas on $\HH_{ab}(H)$. In this paper, we will 
identify freely the monomial $x^p y^q$ with the couple $(p,q)$ 
and therefore the expression ``staircase of monomials'' will make sense.
More generally, we will transpose unscrupulously the definitions
between couples of integers and monomials. If $E$ is a staircase,
then the vector space  $I ^E$ generated by the monomials which are not 
in $E$ is an ideal and reciprocally, every monomial ideal is an 
ideal $I^E$ for a unique staircase $E$. The subscheme  $Z(E)$
whose ideal is $I ^E$ is in $\HH_{ab}(H)$ if and only if 
$E$ has $h_i$ elements in degree $i$.

\subsection{The possible definitions}
\label{defposs}
In this section, we explain that there
are three ways to describe the stratification on $\HH_{ab}(H)$. The 
fact that the  definitions coincide is shown in [Y].
We call a cell of this stratification a ``Schubert cell'' 
on $\HH_{ab}(H)$ as it is a Schubert cell when $\HH_{ab}(H)$ is 
a Grassmannian (section \ref{ssgrass}). Moreover each cell contains
a unique subscheme $Z(E)$ and we will denote this cell by $C(E)$.
\nl
{\bf The Grobner basis point of view.} 
The theory  
of Gr\"obner basis associates with every  
ideal in $k[x,y]$ a monomial ideal (with respect  
to the monomial order chosen in section \ref{ssnot}) called initial ideal  
and $C(E)$ is the locus in $\HH_{ab}(H)$ parametrizing the  
ideals whose initial ideal is $I^E$. For the reader not familiar with 
Gr\"obner basis, we give a characterization which is sufficient for the  
sequel. Let $m_1, m_2, \dots$ be the monomials which don't belong to $E$.  
An ideal of $\HH_{ab}(H)$  
is in $C(E)$ if, regarding it as a $k$-vector space, it admits a base  
$f_1, f_2, \dots$ where $f_i=m_i+R_i$, $R_i$ being a linear combination  
of monomials strictly smaller than $m_i$.  

\brq 
When the product $a.b$ is zero,   
there is at most one staircase compatible  
with the Hilbert function $H$ (i.e. such that $E$ has $h_i$ elements in 
degree $i$) and $\HH_{ab}(H)$ is either empty  
or is reduced to a unique cell (it is non empty when $H$ is  
a decreasing sequence).  These  cases are not 
relevant for  the incidence problem and we suppose from now on $ab\neq 0$  
and $a>0$.  
\erq

\noindent
{\bf The Schubert cells point of view.}
Let  $H=(h_0,\dots,h_r,0,0,\dots)$,
and  $Z \in \HH_{ab}(H)$. The ideal $I_{n}(Z)$ being a vector space 
of codimension $h_n$ in $k[x,y]_n$, it corresponds to a point
$p_n$ in a Grassmannian $G_n$.  
So $I(Z)=\oplus I_n(Z)$ corresponds to a point $p=(p_0,\dots,p_r)$ 
in the product of
Grassmannianns $G_0\x G_1 \dots \x G_r$.
This set theoretical representation turns out to be a closed
embedding $g:\HH_{ab}(H)\fd G_0\x G_1 \dots \x G_r $. 
Let $V_i$ be the 
subspace generated by the $i$ smallest monomials
of $k[x,y]_n$. The flag $V_0 \inc V_1 \dots  \inc V_{n+1}$ defines 
by a classical construction 
Schubert cells on $G_n$ which stratify it. 
The stratification we consider on $\HH_{ab}(H)$ is the stratification
whose strata are the locally closed subschemes 
$g^{-1}(C_0\x C_1\x \dots \x C_r)$, where  $C_i$ is a Schubert cell
on $G_i$. 
\nl
{\bf The Bialynicki-Birula point of view.}
Let $X$ be a smooth projective variety 
over $k$
admitting an action of 
the torus $k^*$. Suppose that the action has a 
finite number of fixed points $x_1, \dots,x_n$. Let 
$T_{X,x_i}^+$ be the part of the tangent space to $x_i$ in 
$X$ where the weights of the $k^*$-action are positive, and 
let $X_i:=\{x \in X, \lim_{t \fd 0}(t.x)=x_i\}$. Then 
a theorem of Bialynicki-Birula asserts that the $X_i$ are a cellular 
decomposition of $X$ in affine spaces and satisfy $T_{X_i,x_i}=T_{X,x_i}^+$.
In our case, fixing two integers $p$ and $q$ with $ap+bq>0$, 
the torus $k^*$ acts on $k[x,y]$ by $t.x=t^p.x$ and $t.y=t^{q}y$. 
This action induces an action of $k^*$ on $\HH_{ab}(H)$. 
The fixed points of $\HH_{ab}(H)$ under $k^*$ are 
the monomial subschemes $Z(E)$.
Applying the Bialynicki-Birula theorem to the action of $k^*$ on $\HH_{ab}(H)$,
we get a stratification and $C(E)$ is the cell associated with the
fixed  point $Z(E)$.

\subsection{Incidence relations}
In this section, we recall the main theorem of [Y]
about incidence relations.
\\
The monomials 
of $k[x,y]$ can be ordered in 
an infinite sequence $m_0<m_1<\dots$ thanks to the monomial order. 
Let $S_E$ be the function
from $\NN$ to $\NN$ defined by $S_E(k)=$number of monomials in  
$E$ smaller or equal to $m_k$.
\bthm
\label{resultJ}
If $\overline{C(E)}\cap C(F)\neq \emptyset$ then $S_E\geq S_F$.
\end{thm}

\brq
In [Y], the result says $\leq $ instead of $\geq$ as above. 
The reason is that we
have slightly changed the definition of $S_E$ to get a shorter
presentation. Moreover, 
the paper deals with  the homogeneous 
case $a=1,b=-1$ but the extension is
immediate.
\erq

\noindent
The condition is not sufficient: if $E$ and $F$ are the 
staircases whose monomial ideals are 
$I^E=(y^4,xy^2,x^2y,x^5)$ and $I^F=(y^5,xy^2,x^3)$ 
then $S_E \geq S_F$ but the incidence $\overline{C(E)}\cap C(F)\neq 
\emptyset$ 
is not fulfilled([Y]).
 
\section{A necessary condition for weak incidence}
\label{s2}

\subsection{Statement of the result}
In this section, we give a necessary condition on two staircases 
$E$ and $F$ to fulfill the relation $\overline{C(E)}\cap C(F)
\neq \emptyset$. The condition will rely on the notions 
of system of arrows  and of dual staircase that we introduce now. 
\bdefi
An arrow on $E$ is a couple of monomials $(P,Q)$, ($P$ being the origin
of the arrow and $Q$ being the end of the arrow), such that $P$ is  
in $E$ and  such that the vector $\overrightarrow{PQ}$  of $\NN^2$ 
is a negative multiple of $(a,b)$ (recall that we have adopted 
the convention $a>0$). If $f=(P,Q)$, the multiplication by a monomial 
$m$ of $f$ is the arrow $m.f:=(mP,mQ)$.
An arrow 
$f=(P,Q)$ is shorter that  $f'=(P',Q')$ if $\overrightarrow{PQ}=\lambda
(a,b)$, $\overrightarrow{P'Q'}=\lambda'
(a,b)$, and if the absolute values of $\lambda$ and $\lambda'$ 
verify $|\lambda| \leq |\lambda'|$.\\
 A system of
arrows on $E$ is the data for each element 
 $p \in E$ of an arrow $(p,f_p)$
on $E$ such that
 set $S$ of arrows 
parametrized by $E$ satisfy the
following   conditions: 
\bit
\item
$p \neq q \Rightarrow f_p \neq f_q$
\item
the arrows are compatible with division, which means that   
$\forall f=(P,Q)\in S$, $\forall Q'$ monomial of 
$k[x,y]$ dividing $Q$,  $\exists g \in S$ 
such that $g$ has end $Q'$ and such that $g$ is shorter 
than $f$. 
\eit
\edefi

With a staircase $E$ and a system of arrows $S$, one can define a new 
staircase $S(E)$ as follows.
\bprop \label{S(E)est_un_esc}
Let $E$ be a staircase and $S$ a system of arrows on $E$. Let $S(E)$
be the set 
of monomials which are the end of an arrow. Then $S(E)$ is a staircase.
\eprop
We have to verify that if $m$ and $Q'$ are two monomials with
$mQ' \in S(E)$ then $Q' \in S(E)$. By definition of 
$S(E)$, $Q:=mQ'$ is the end of an arrow. 
The compatibility of $S$ with
division shows that $Q'\in S(E)$.
\findem

\bdefi
If $E$ and $F$ are two staircases, and if $S$ is a system of arrows
on $E$ such that $F=S(E)$, we will say that $S$ is a system of arrows
from $E$ to $F$. 
\edefi

\bdefi
The box of size $M\x N$ is $B_{M\x N}:=\{ (x,y) \in \NN^2, x<M, y<N\}$. 
If a staircase $E$ is included in $B_{M\x N}$, the dual of $E$ 
in the box is by definition the
set $E^{\nu}_{MN}:=\{(x,y)\in \NN^2 such\ that\ \forall (e_1,e_2)
\in E, e_1+x<M-1 \mbox { and } e_2+y<N-1\}$.
In the sequel, we will often write $E^{\nu}$ instead of 
$E^{\nu}_{MN}$
for simplicity. 
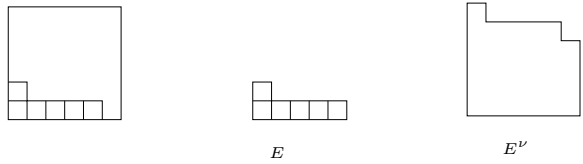
\edefi

  \begin{figure}[h] 
     \begin{center}
        \input{escetsondual.pstex_t}
     \end{center} 
  \end{figure}

\bthm 
\label{thm} 
Let $E$ and $F$ be two staircases 
included in a box $B_{M\x N}$. Let $E^{\nu}$ and  
$F^{\nu}$ be their respective dual in $B_{M\x N}$. 
If the incidence $\overline{C(E)}\cap C(F)\neq  
\emptyset$  
is fulfilled, then  
\ben 
\item 
there exists a system of arrows $S$ on $E$ such that $S(E)=F$ 
\item 
there exists a system of arrows $S^\nu$ on $E^\nu$ such that 
$S^\nu(E^{\nu})=F^{\nu}$ 
\een 
\ethm 
see section \ref{proof} 
\findem

\brq 
One can prove that the existence of 
$S^{\nu}$ does not depend on the choice of the box $B_{M\x N}$
in which we construct the dual staircase $E^{\nu}$.
\erq

\noindent
The following example shows that conditions 1) and 2) are not equivalent.
\bex
\label{condnoneq}
Let $E$ and $F$ be the staircases whose monomial ideals are   
$I^E=(y^4,xy^2,x^2y,x^5)$ and $I^F=(y^5,xy^2,x^3)$. Let $(a,b)=(1,-1)$.
There exist 
a system $S$ on $E$ such that $S(E)=F$ but there is no system
of arrows $S$ on $E^{\nu}_{5\ 5}$ such that $S^{\nu}(E^{\nu}_{5\ 5})=F^{\nu}_{5\ 5}$.
\eex
The set $S:=\{(x^3,x^2y),(x^4,y^4),(p,p), p\in E, p\neq x^3, p\neq x^4\}$
is a system of arrows on $E$ such that $S(E)=F$. 
The dual staircases $E^{\nu}$ and $F^{\nu}$ in $B_{5\x 5}$ are such that 
$I^{E^{\nu}}=(y^4,x^3y^3, x^4y,x^5)$ and $I^{F^{\nu}}=
(y^5,x^2y^3, x^4)$.
Suppose that there exists  
a system of arrows $S^{\nu}$ on $E^{\nu}$ such that
$S^{\nu}(E^{\nu})=F^{\nu}$.  
Because $x^2y^3$ is in $E^{\nu}$ but not in $F^{\nu}$, $S^\nu$
must contain   an arrow $f$ from $x^2y^3$ to a point $p> x^2y^3$,
$p\in F^\nu$ ie.
$p=xy^4$.  The monomial  $x^3y^2 \in F^{\nu}$ so there is an arrow $f'$ from  
a monomial $p \in E^{\nu}$ and $\leq x^3y^2$ to $x^3y^2$, i.e. $p=x^3y^2$. 
The compatibility with division and the arrow  
$f'$ show that all 
monomials dividing $x^3y^2$ admit an arrow to themselves.  
The compatibility with division by $y$ applied to the  
arrow $f$ ensures the existence of an arrow from $xy^3$ to itself.  
It follows that the compatibility with division by $x$  
applied to $f$ cannot be satisfied.  
\findem

\subsection{Proof of theorem \ref{thm}}
\label{proof}
\subsubsection{Proof of point 2}
To show that there exists a system 
$S^\nu$ on $E^\nu$ such that
$S^\nu(E^{\nu})=F^{\nu}$, 
we will exhibit a system $S^c$ on the complementary $E^c$
of $E$ in $\NN^2$ such that $S^c(E^c)=F^c$. 
So the first step is to explain what is a system on $E^c$, 
and to reduce the proof to a construction of such a system $S^c$. 
The reduction is given in proposition \ref{reduction}

\begin{defprop}
An arrow on $E^c$ is a couple of monomials $(P,Q)$ 
such that $P \notin E$ and  such that the vector $\overrightarrow{PQ}$  of
$\NN^2$   is a positive multiple of $(a,b)$.
A system of arrows on $E^c$ is the data for each element 
$p \notin E$ of an arrow $(p,f_p)$ on $E^c$ such that
the
set $S^c$ of arrows 
parametrized by $E^c$ satisfy the
following   conditions: 
\bit
\item
$p \neq q \Rightarrow f_p \neq f_q$
\item
the arrows are compatible with multiplication, which means that   
$\forall f=(P,Q)\in S^c$, $\forall m$ monomial of 
$k[x,y]$,  $\exists g \in S^c$ 
such that $g$ has end the product $mQ$ and
such that $g$ is shorter than $f$.
\eit
If $S^c$ 
is a system of arrows on $E^c$, then the set $$S^c(E^c)
:=\{(x,y) \mbox{ which are the end of an arrow of } S^c\}$$ is 
the complementary of a staircase $F$. The system $S^c$ will be 
called a system of arrows from $E^c$ to $F^c$.
\end{defprop}
\noindent
The proof is similar to the proof of proposition \ref{S(E)est_un_esc}
\findem
\nl
The next lemma identifies systems on $E^c$ and systems on $E^{\nu}$.
Let $\goth S^{\nu}$ (resp. $\goth S^c$) 
be the set of systems of arrows on $E^{\nu}$ (resp. on $E^c$). 
Let $\phi:B_{M\x N} \fd B_{M\x N}$, $(x,y)\mapsto (M-1-x,N-1-y)$
be the dualizing map. For an arrow $f=(p,q)$  on
$E^{\nu}$, if $ q \in B_{M\x N}$, 
we define $\phi(f):=(\phi(p),\phi(q))$. Let 
$$\goth S^{\nu}_{MN}:=\{ S^\nu \mbox{ system of arrows on $E^\nu$ s.t.  } 
s=(p,q) \in S^{\nu} \Rightarrow q \in B_{M\x N}\}.$$
If  $S^{\nu} \in \goth S^{\nu}_{MN}$,  
let $\psi(S^{\nu}):=\{\phi(f), f\in S^{\nu}\} \cup \{(p,p), \ p\in
E^c\setminus \phi(E^{\nu})\}$. 
 
\blm
\label{descimage}
The map $\psi: \goth S^{\nu}_{MN} \fd \goth S^c$  is well defined, 
(ie. $\psi (S^{\nu})$ is a system of arrows on $E^c$) and injective. 
The image $Im \ \psi$ contains  the systems $S^c$ such
that:  for all monomial $m \notin B_{M\x N}$, the arrow of 
$S^c$ whose origin is $m$ is the arrow $(m,m)$.  
Moreover if $S^{\nu}$ and $S^c$ are two systems
such that  $\psi(S^{\nu})=S^c$, then
$S^{\nu}(E^{\nu})=(S^c(E^c)^c)^{\nu}$. 
\elm
it consists in a sequence of easy 
combinatorial verifications 
which are left to the reader.
\findem
\nl
Here is a another criterion to check whether a system $S^c$ 
on $E^c$ is in $Im\ \psi$.

\blm
Let $E$ be a staircase included in 
the box $B_{M\x N}$, $S^c$ be a system 
of arrows on $E^c$, and $F$ be the staircase such that
$S^{c}(E^c)=F^c$. If $F \inc B_{M\x N}$, then $S^c \in Im\ \psi$. 
\elm
Let $p_0, \dots ,p_s$ be the set of monomials 
of $k[x,y]_s$. Up to reordering, one can suppose 
that $p_0>p_1>\dots>p_s$. There exist integers 
$l$ and $m$ such that $k[x,y]_s
\setminus B_{M\x N}=\{p_0,\dots, p_l\} \cup \{p_m,\dots, p_s\}$,
$l<m$, $-1\leq l \leq s$, $0\leq m\leq s+1$. 
By \ref{descimage}, to prove the proposition, we have to show 
that $\forall i \in \{0,\dots, l-1,l,m,m+1,\dots, s\}$, 
the arrow $(p_i,p_i)$ is in $S^c$. 
If $m\leq s$,  $p_s \notin B_{M\x N}$ so 
$p_s \in E^c$. It follows that  $S^c$ 
contains an arrow $(p_s,*)$ by definition of a system 
on $E^c$ and necessarily 
$*=p_s$ since $*\leq p_s$. Now, if $m<s$, $p_{s-1}
\in E^c$, so there is in $S^c$ an arrow $(p_{s-1},*)$.
But distinct arrows have distinct ends so 
$*=p_{s-1}$. By decreasing 
induction, one shows that $S^c$ contains $(p_m,p_m),
(p_{m+1},p_{m+1}),
\dots, (p_s,p_s)$. \\
Similarly, if $p_0 \notin B_{M\x N}$ then $p_0 \in F^c$.
By definition of $S^c(E^c)$, this means that $S^c$ contains an
arrow $(*,p_0)$  and 
the only possibility is $(p_0,p_0)$. If the monomial $p_1
\notin B_{M\x N}$, then $p_1$ 
is the end of an arrow $(*,p_1)$ and we must 
have $(*=p_1)$ because the arrow with origin 
$p_0$ is already determined.  By induction, $S^c$ contains 
the arrows $(p_0,p_0), \dots, (p_l,p_l)$. 
\findem

\bprop
\label{reduction}
To prove the 
existence of a 
system $S^{\nu}$ such that $S^{\nu}(E^{\nu})=F^{\nu}$,
it suffices to exhibit a system $S^c$ on $E^c$ 
such that $S^c(E^c)=F^c$.
\eprop
it is an immediate consequence of the last two lemmas. 
\findem
\nl
The goal of the 
next proposition is to exhibit such a system $S^c$.  
The cell $C(E)$ is an affine space as a cell 
of a Bialynicki-Birula stratification. So, if $\overline{C(E)}$ meets 
$C(F)$, there exists a morphism 
$f$ from the affine line $\AA^1=\s k[t]$ to 
$\overline{C(E)}$ such that the image of the generic point is 
in $C(E)$ and such that the limit 
of $f(t)$ when $t$ tends to infinity is $X_{\infty} \in C(F)$. 
Let $p_1,\dots,p_s$ be the points in $\AA^1$ whose image by 
$f$ is in $\overline{C(E)}-C(E)$. By the universal property of $C(E)$,
the morphism $f$ is defined by a closed subscheme $U$ of 
$\AA^1-\{p_1,\dots,p_s\}\x \s k[x,y]$. Let 
$I(t)\inc k[x,y][t]$ be the ideal defining the closure of 
$U$ in $\AA^1\x \s k[x,y]$. 
 The staircase
over the generic point being $E$, for 
 each monomial $m \notin E$, there
exists a quasi-homogeneous element 
 $P \in I(t)$ 
with initial monomial $m$. 
Among all possible $P$, choose one as follows. 
Let $d_t(P)$ be the degree of $P$ in $t$, and $d_{in,t}(P)$ 
be the degree in $t$ of the initial coefficient  of $P$. Let $S_m$ 
be the set of $P$ such that $\Delta(P):=d_t(P)-d_{in,t}(P)$ is minimal. 
For a fixed $P$, denote by  $val(P)$ the greatest monomial of $P$  
whose coefficient has degree $d_t(P)$. Now choose a $P$ in 
$S_m$ such that $val(P)$ is minimum. 
Call this element $P(m)$ and put $e(m):=val(P(m))$.

\bprop \label{Sc_est_un_syst_de fleches}
Let $S^c:=\{(m,e(m)),\ m \in E^c\}$. Then $S^c$ is a system 
of arrows on $E^c$ and $S^c(E^c)=F^c$. 
\end{prop}
To prove the proposition, we need two lemmas.

\blm \label{lemme_technique}
Let $P$ be a quasi-homogeneous element of $I(t)$ with $in(P)=m$ and 
$val(P)=n$. Let $Q$ be a quasi-homogeneous 
element of $I(t)$ such that $in(Q)<m$  
and $val(Q)=n$. Let $b(t)$ be the coefficient
of $n$ in $P$, $d(t)$ be the coefficients of 
$n$ in  $Q$. Then $R=d(t)P-b(t)Q$ is an element of $I(t)$
with initial term $m$ satisfying 
$\Delta(R)<\Delta(P)$ or ($\Delta(R)=\Delta(P)$ and $val(R)<val(P)$).
\elm
the fact that $R$ is in $I(t)$ with initial term $m$ is obvious
and it remains to calculate $\Delta(R)$. 
By construction, $d_t(R)\leq d_t(P)+d_t(Q)$ and 
$d_{in,t}(R)=d_{in,t}(P)+d_t(Q)=d_t(P)-\Delta(P)+d_t(Q)$ 
so $\Delta(R)\leq \Delta(P)$. If $\Delta(R) = \Delta(P)$, 
then $val(R)$ is the greater monomial of $R$ whose coefficient 
in $t$ has degree $d_t(P)+d_t(Q)$.
By construction, $R$ has no term in $n$ and 
all monomials greater than $n$ have degree in $t$  strictly less 
than $d_t(P)+d_t(Q)$. So, the greatest 
monomial having degree $d_t(P)+d_t(Q)$ is smaller than $n$, showing
$val(R)<val(P)$. \findem 

\blm
1) The map $e:E^c\fd \NN^2$ sending $m$ to $e(m)$ is injective. \\
2) $Im\ e=F^c$.
\\
3) The inverse map to $e:E^c \fd F^c$ sends a monomial $n\notin F$ to the
smallest monomial 
 $m \notin E$ for which there exists a quasi-homogeneous
polynomial $P \in I(t)$ with $in(P)=m$
 and $val(P)=n$. 
\elm
1) If $e$ was not injective, there would exist two monomials $m_1$ and 
$m_2$ 
with $m_1>m_2$ and $val(P(m_1))=val(P(m_2))$. 
Then the preceding lemma would construct 
from $P(m_1)$ and $P(m_2)$ an element $R\in I(t)$ 
contradicting the minimality of $P(m_1)$. \\
2) By definition of the flat limit, for any $m\notin E$, 
the limit $P_{\infty}(m)$ of $P(m)/t^{d_t(P(m))}$ as  $t \rightarrow 
\infty$ is an element of $I(X_{\infty})$. It follows  
by definition of $X_{\infty}$ that the initial 
monomial $e(m)$ of $P_{\infty}(m)$  is not in $F$.
Conversely, we show  that an element $n \notin F$ 
can be written $e(m)$ for some $m \in E^c$. 
Let $S$ (resp. $T$) be the set of elements of $E^c$ 
(resp. of $F^c$) with the same degree as $n$. The flatness
shows that $S$ and $T$ have the same (finite) cardinal. 
By the above, $e$ is 
an injection from $S$ into $T$, so is a bijection. In particular, $n
\in T$ can be written $e(m)$ for some $m$. \\
3) Let $n=e(m)$ be a monomial in $F^c$. The polynomial $P(m)\in I(t)$
verifies $in(P(m))=m$ and $val(P(m))=n$. Suppose that there exists 
$m'<m$ and a quasi-homogeneous polynomial $P'\in I(t)$ such that $in(P')=m'$
and $val(P')=n$. 
 Then lemma \ref{lemme_technique} would construct from $P(m)$ and $P'$ 
a polynomial $R$ contradicting the minimality of $P(m)$. 
\findem
\nl
Proof of proposition \ref{Sc_est_un_syst_de fleches}.
Point 1) of the last lemma shows that each point is the end of at most
one arrow. 
To prove that the set is a system of arrows, 
the compatibility
with multiplication by a monomial $m'$ remains to be shown. 
 Let $f=(m,e(m))$ be an arrow. Then
$P(m)\in I(t)$ verifies $in(P(m))=m$
 and $val(P(m))=e(m)$. The polynomial
$m'P(m)\in I(t)$ 
verifies $in(m'P(m))=m'm,val(m'P(m))=m'e(m)$. Then
point 3)
 shows that the arrow $g$ ending at $m'e(m)$ has an origin smaller or
equal to $xm$, i.e. $g$ is shorter than $f$. 
So $S^c$ 
is a system of arrows on $E$, and $S^c(E^c)=F^c$ by 2). 
\findem

\subsubsection{Proof of point 1}
It will essentially rely on the proof of 
point 2 thanks to a notion of duality
corresponding to the notion of linkage in 
a complete intersection, introduced in [PS].

\bdefi
Let $Z$ be a zero-dimensional scheme included 
in the complete intersection $Y=(x^M,y^N)$. 
The scheme $Z^\nu$ defined by 
the ideal $(I(Y):I(Z))$ is called the 
scheme obtained by linkage in the complete intersection
$(x^M,y^N)$. Let $\HH(Y)$ be the reduced subscheme of $\HH$ 
parametrizing the subschemes included in $Y$. 
\edefi

\bprop \label{dualite}
The morphism $\phi:\HH(Y) \fd \HH(Y)$, $Z \mapsto Z^{\nu}$ is well defined
and sends $C(E)$ to $C(E^{\nu}_{MN})$.  
Moreover $Z^{\nu\nu}=Z$. 
\eprop
if $Y$ is a $0$-dimensional Gorenstein scheme, and if $Z\inc Y$,
then the ideal $(I(Y):I(Z))$ defines a scheme of degree 
$deg(Y)-deg(Z)$. A complete
intersection being Gorenstein, the degree of 
$Z^\nu$ is constant on the connected components of $\HH(Y)$ and 
the morphism is well defined. 
\\
Let $Z \in C(E)$, let $in(I(Z^{\nu}))$ be the set of initial
monomials of the elements of $I(Z^{\nu})$, and let $e=x^{e_1}y^{e_2}
\in in(I(Z^{\nu}))$. The inclusion 
$in(I(Z^{\nu}))in(I(Z)) \inc in(I(Z)I(Z^{\nu})) \inc in(x^M, y^N)$
shows that: $\forall (\alpha, \beta)\notin E$, $e_1+\alpha \geq M$ or 
$e_2+\beta\geq N$. In other words, $in(I(Z^{\nu}))\inc (E^{\nu}_{MN})^c$. 
The complementary sets $in(I(Z^{\nu}))^c$ and $E^{\nu}_{MN}$ have the 
same cardinal $MN-card(E)$ so $in(I(Z^{\nu}))=(E^{\nu}_{MN})^c$. 
This shows that 
$\phi$ sends $C(E)$ to $C(E^{\nu}_{MN})$. 
\\
Finally, we have by definition
of the dual, $Z^{\nu\nu}\inc Z$. But $Z$ and $Z^{\nu\nu}$ have the 
same degree so they are equal. 
\findem

\bcor
Let $E$ and $F$ be two staircases of the same cardinal $n$. 
Then $\overline {C(E)} \cap C(F) \neq \emptyset \Rightarrow
\overline {C(E^{\nu}_{nn})} \cap C(F^{\nu}_{nn}) \neq \emptyset$.
\ecor
\textit{Proof}: apply the dualizing morphism $\phi$ of the last proposition with 
$M=N=n$. 
\findem
\nl
Now, we come to the proof of point 1 of theorem \ref{thm}. 
If $\overline {C(E)} \cap
C(F)
 \neq \emptyset$, then $\overline {C(E^{\nu}_{nn})} \cap C(F^{\nu}_{nn}) \neq
\emptyset$.
 Then, point 2 states that there exists a system of arrows $S$ on
$(E^{\nu}_{nn})^{\nu}_{nn}=E$ such that $S(E)=(F^{\nu}_{nn})^{\nu}_{nn}=F$.
\findem 

\section{Comparison to the known results}
\subsection{Grassmannians}
\label{ssgrass}
In the case  $H=(1,2,3,\dots,k-1,k,k+1-n,0,0,\dots,0)$,
$\HH_{1,-1}(H)$ is isomorphic to the Grassmannian $G(n,k+1)=G(n,k[x,y]_{k})$.
The isomorphism consists in associating  $I=\oplus I_d$ with the vector 
space $I_{k}$ and the inverse isomorphism associates $I_k$ with the ideal 
$I=I_k\oplus k[x,y]_{k+1}\oplus k[x,y]_{k+2} \oplus \dots$ . A
non increasing 
sequence of  non negative integers $(p_1,\dots,p_n)$ with
$p_1 \leq k+1-n$ defines a staircase  $E(p_1,\dots,p_n)$ which contains the
monomials of degree at most $k-1$  and the monomials of degree $k$ 
which are not in
$\{x^{n-i+p_i}y^{k-n+i-p_i}$, $1\leq i\leq n\}$. 
Let  $V_i$ be  the vector space generated by the
$i$ smallest monomials of degree $k$.
The  Schubert
cell $C(p_1,\dots,p_n)$ on $G(n,k[x,y]_{k})$ is by definition 
$C(p_1,\dots,p_n):=\{ W \in G(n,k[x,y]_{k}), dim(W \cap V_j)=i\mbox { if }
k+1-n+i-p_i\leq j < k+1-n+i-p_{i+1} \}$.
By construction, the Schubert
 cells $C(E(p_1,\dots,p_n))$ and
$C(p_1,\dots,p_n)$ correspond  under the above
 isomorphism. A classical
result on Grassmannians asserts  that the closure of
 a cell
$C(p_1,\dots,p_n)$ is the union of the cells  $C(q_1,\dots,q_n)$ for
 which
$q_i\geq p_i$ for all $i$.  The following result says that in the
 particular
case of Grassmannians, the two points of our  
necessary condition coincide, and they coincide  with the
known  incidence condition on Grassmannians.

\bprop
\label{idemG}
The four following conditions are equivalent:
\ben
\item
$\overline {C(p_1,\dots, p_n)} \cap C(q_1,\dots, q_n) \neq \emptyset$
\item
$q_i\geq p_i$ for all $i$.
\item 
There exists a  system $S$ on $E(p_1,\dots,p_n)$ such that 
$S(E(p_1,\dots,p_n))=E(q_1,\dots,q_n)$
\item 
There exists a  system $S^{\nu}$ on $E(p_1,\dots,p_n)^{\nu}$ such that  
$S^{\nu}(E(p_1,\dots,p_n)^{\nu})=E(q_1,\dots,q_n)^{\nu}$ \een
\eprop
$1 \Leftrightarrow 2$ is  a classical result on Grassmanians.\\
$1 \Rightarrow 3$ 
and $1 \Rightarrow 4$ is exactly theorem \ref{thm}.
\\
$4 \Rightarrow 1$. If 4) is true, we have a system $S^{\nu}$ 
on $E^{\nu}$, which can be identified with a system $S^c$ on 
$E^c$ (thanks to lemma \ref{descimage}) such that $S^c(E(p_1,\dots,p_n)^c)
=E(q_1, \dots, q_n)^c$. 
For each arrow $f$ of $S^c$ in degree $k$ with origin 
$o(f)$ and end $g(f)$, let $e(f):=o(f)+tg(f)\in k[x,y,t]$ and denote by $V_k$ 
the $k[t]$-module generated by $\{e(f),\ f\in S^c \mbox { of degree } k\}$. 
Let $I(t)\inc k[x,y][t]$ be the ideal defined 
by $I(t):=V_k\oplus k[x,y]_{k+1}[t] \oplus k[x,y]_{k+2}[t] \oplus \dots$. The
scheme whose ideal is $\lim_{t \fd
 \infty}I(t)=I^{E(q_1,\dots,q_n)}$ is  in
$\overline{C(E(p_1,\dots,p_n))} \cap
 C(E(q_1,\dots,q_n))$. \\
$3 \Rightarrow 4 $. We have $4 \Rightarrow 3$. Using the duality of proposition
\ref{dualite} to exchange the roles of $E(p_1,\dots,p_n)$
and $E(p_1,\dots,p_n)^{\nu}$, we have $3 \Rightarrow 4$.
\findem

\subsection{Comparison to Yameogo's result} 
The next result asserts that 
any of the two conditions of theorem
\ref{thm} implies strictly the condition of theorem \ref{resultJ}.

\bprop
\label{amelJ}
Let $S$ be a system of arrows on $E$. Then $S_{S(E)}\leq S_E$. 
Moreover, there exist staircases $E$ and $F$ such that:
\bit
\item
$S_E \geq S_F$
\item
there is no 
system of arrows $S$ on $E$ verifying $F=S(E)$
\eit
Let $S^{\nu}$ be a system of arrows on $E^{\nu}$. Then
$S_{S^{\nu}(E^{\nu}_{MN})^{\nu}_{MN}}\leq S_E$.  
 Moreover, there exist staircases
$E$ and $F$ 
such that:
\bit
\item
$S_E \geq S_F$
\item
there is no 
system of arrows $S^{\nu}$ on $E^{\nu}_{MN}$ 
verifying $F^{\nu}_{MN}=S^{\nu}(E_{MN}^{\nu})$
\eit
\eprop
The passage from $E$ to $S_E$ can be described dynamically as 
follows. First, you  suppress all the arrows which go from 
a point $P$ to itself. Then, if there is an arrow $f_1$ from 
a monomial $P$ to a monomial $Q$ 
and an arrow  $f_2$ from $Q$ to a monomial $R$, you suppress 
the arrows $f_1$ and $f_2$ and you replace them by the arrow 
from $P$ to $R$. After a finite number of such operations,
you come to a finite set of arrows $g_1, \dots, g_s$ 
with distinct ends
such that the origins of the arrows 
belong to $E$ and such that the ends of the arrows 
do not belong to $E$. Now suppress in $E=:E_1$ the origin of 
$g_1$ and replace it by the end  of $g_1$. You obtain in this way 
a subset $E_2$ of $\NN^2$. Replacing in $E_2$ the origin of 
$g_2$ by the end of $g_2$, a subset $E_3$ of $\NN^2$ is obtained. 
Continuing in this way, you finally reach the subset $E_{s+1}=S(E)$. 
We have by construction $S_{S(E)}=S_{E_{s+1}}\leq S_{E_{s}} \leq \dots
\leq S_{E_2} \leq S_{E_1}=S_E$.
\\
The same reasoning shows that the inequality 
$S_{S^{\nu}(E^{\nu}_{MN})^{\nu}_{MN}}\leq S_E$ always holds. 
\\
For the second part of the proposition, we 
consider the staircases 
$E$ and $F$ 
whose monomial ideals are  
$I^E=(y^4,xy^2,x^2y,x^5)$ and $I^F=(y^5,xy^2,x^3)$.
They verify $S_E\geq S_F$ but example \ref{condnoneq} shows that 
there is no  
system of arrows $S^{\nu}$ on $E^{\nu}_{5\ 5}$ verifying $F^{\nu}_{5\ 5}=
S^{\nu}(E^{\nu}_{5\ 5})$.
The dual example $I^E=(y^4,x^3y^3, x^4y,x^5)$ and $I^{F}=
(y^5,x^2y^3,x^4)$ verifies $S_E\geq S_F$, but 
there is no   
system $S$ on $E$ such that $F=S(E)$.
\findem

\section{An example} \label{sex}
For a small Hilbert function, it is easy to 
show that the necessary condition given by our theorem is 
in fact sufficient. It follows that we can solve  the 
weak incidence problem on $\HH_{ab}(H)$ for small $H$. 
In this section, we illustrate this fact in the case 
$(a,b)=(1,-1)$ and 
$H=(1,2,3,2,1,0,0,\dots)$.
\nl
To solve the problem, we will use 
sufficient incidence conditions given
by $[Y2]$. We could deal 
without these conditions, but we use them to shorten the proofs.
\\
There are nine possible staircases corresponding to the Hilbert
function $H$, all of them drawn and named in the following figure.  

  \begin{figure}[h] 
     \begin{center}
        \input{stratesG12321l.pstex_t}
     \end{center} 
  \end{figure} In the figure, a square whose position is $(i,j)$ 
relative to the square at the left bottom represents the monomial 
$x^iy^j$.
We know that $\overline {C(E_1)} \cap C(E_2)\neq \emptyset 
\Rightarrow S_{E_1}\geq S_{E_2}$. It follows that the only 
possible relations are $$\{(E_{gen},*);(a,c);(a,d);(a,e);(a,f); (a,g);
(a,h);(b,d);(b,e);(b,f);(b,g);$$ $$(b,h);(c,e);(c,f);(c,g);(c,h);
(d,f);(d,g);(d,h);(e,g);(e,h); (f,h);(g,h)\}$$
where $*$ denotes any staircase and where  we 
have adopted the notation $(E_1,E_2)$
for the incidence relation $\overline {C(E_1)} \cap C(E_2)\neq
\emptyset$.\\
In the figure, a thick arrow between two staircases $E_1$ and $E_2$
means that the relation $\overline {C(E_1)} \supset C(E_2)$ occurs and 
these strong incidence relations come from [Y2].   It follows that 
if two staircases $E_1$ and $E_2$ are linked by a sequence of 
thick arrows, then the relation $(E_1,E_2)$ holds. The cases 
then remaining from the above list are:
$$
\{ (a,e);(b,d);(b,f);(c,e);(c,g);(d,f)\}
$$
Example \ref{condnoneq} shows that  $(c,g)$ does not hold. The five
remaining cases are compatible with our necessary
condition. In particular, there are  
systems of arrows on the dual staircases $a^{\nu},
b^{\nu},b^{\nu},c^{\nu},d^{\nu}$
coming from 
the condition 2) of theorem \ref{thm}. We interpret 
these systems as systems on  $a^{c},
b^{c},b^{c},c^{c},d^{c}$
thanks to lemma
\ref{descimage}. The following
figure shows the systems obtained in this way. 

  \begin{figure}[h] 
     \begin{center}
        \input{systflechesl.pstex_t}
     \end{center} 
  \end{figure}
 To exhibit ideals ``compatible'' with these 
systems, consider the sub-$k[t]$-modules of $k[x,y,t]$
given by their generators: $T_3=<y^3+t^2xy^2+tx^3, x^2y>$,
$T_4=<y^4+t^2xy^3,xy^3+tx^4,x^2y^2,x^3y>$, $P_{\geq
5}=<x^{\alpha}y^{\beta},\ \alpha+\beta\geq 5>$. Consider the ideal
$T:=T_3\oplus T_4\oplus  P_{\geq 5}$ of $k[x,y,t]$. The quotient
$k[x,y,t]/T$  is  flat over
$k[t]$.  The associated universal morphism sends the generic point to a point
in $C(a)$ and can be extended by sending $\infty$ to the point whose 
ideal is the monomial ideal with staircase $e$. This shows that $(a,e)$ 
holds. \\
Similarly, if  
$U_3=<y^3,xy^2+tx^2y+t^2x^3>$, $U_4=<y^4,xy^3,x^2y^2+tx^3y,x^4>$,
$V_3=<y^3+tx^2y,xy^2+tx^3>$,
$V_4=<y^4,xy^3+x^3y,x^2y^2,x^4>$, $W_3=<xy^2,x^2y>$,
$W_4=<y^4+tx^4,xy^3,x^2y^2,x^3y>$, $Z_3=<y^3+txy^2+t^2x^2y,x^3>$,
$Z_4=<y^4,xy^3+tx^2y^2,x^3y,x^4>$,
$U:=U_3\oplus U_4\oplus 
P_{\geq 5}$,
$V:=V_3\oplus V_4\oplus 
P_{\geq 5}$, $W:=W_3\oplus W_4\oplus 
P_{\geq 5}$, $Z:=Z_3\oplus Z_4\oplus P_{\geq 5}$,
the ideals $U,V,W,Z$ show respectively that the relations $(b,d),
(b,f), (c,e), 
(d,f)$ are true. The diagram sums up our results: the relation 
$(E,F)$ holds if and only if the staircases $E$ and $F$ are linked 
by a thin arrow or by a sequence of thick arrows. 
\nl
{\bf Bibliography}:
\nl
[BB]: Bialynicki-Birula, A: Some theorems on actions of 
algebraic groups . \textit{Ann. of Math. 98, 480-497, (1973)}
\\ \noindent
[BB1]: Bialynicki-Birula, A: Some properties of the decompostions 
of algebraic varieties determined by actions of a torus. \textit{
Bulletin de l'acad\'emie Polonaise des sciences, Serie des Sciences
math. astr. et phys. 24, \#9, 667-674, (1976)}\\ \noindent
[Br]: Brion, M: Equivariant Chow groups for torus actions, \textit{ 
Transformations Groups 2, 225--267, (1997)}\\ \noindent 
[ES]: Ellingsrud, G. and Stromme, S.: On the homology of the 
Hilbert scheme of points in the plane, \textit{Invent. Math. 87 (1987),
 no. 2, 343--352}  \\ \noindent
[GH]: Griffiths, P. and  Harris, J.: Principles of algebraic geometry, 
\textit{J. Wiley, New York (1977)}\\ \noindent
[Go]: G\"ottsche, L.: Betti numbers for the Hilbert function strata
of the punctual Hilbert scheme in two variables. 
\textit{ Manuscripta Math. 66, 253-259 (1990)}\\ \noindent
[Gr]: Granger, JM.: G\'eom\'etrie des sch\'emas de Hilbert ponctuels,
\textit{ Mem. Soc. Math. Fr, Nouv. ser. 7-12(1982-83)}\\ \noindent
[I]: Iarrobino, A: Punctual Hilbert Scheme. \textit{ Memoirs of 
AMS, vol. 10, \#188, (1977)}\\ \noindent
[M]: Manivel, L.: Fonctions sym\'etriques, polyn\^omes de Schubert
et lieux de d\'eg\'en\'erescence,  \textit{
Cours specialises de la SMF, 3, (1999)}\\ \noindent
[PS]: Peskine C., Szpiro L.: Liaison des vari\'et\'es alg\'ebriques,
\textit{Invent. Math. 26, (1974), 271-302}\\ \noindent   
[Y]: Yameogo, J: D\'ecomposition cellulaire de vari\'et\'es 
 param\'etrant des
id\'eaux homog\`enes de $\CC[[x,y]]$. Incidence 
 des cellules I. \textit{ Compositio Math. 90, \#1, 
 81--98, (1994) }\\ \noindent
[Y2]: Yameogo, J: D\'ecomposition cellulaire de vari\'et\'es 
param\'etrant des id\'eaux homog\`enes de $\CC[[x,y]]$. Incidence 
des cellules II. \textit{ J.reine angew. Math. 450, 123-137,
(1994) }
\end{document}

%% file: escetsondual.pstex_t
\begin{picture}(0,0)%
\includegraphics{escetsondual.pstex}%
\end{picture}%
\setlength{\unitlength}{2072sp}%
\begingroup\makeatletter\ifx\SetFigFont\undefined%
\gdef\SetFigFont#1#2#3#4#5{%
  \reset@font\fontsize{#1}{#2pt}%
  \fontfamily{#3}\fontseries{#4}\fontshape{#5}%
  \selectfont}%
\fi\endgroup%
\begin{picture}(6864,1857)(439,-1411)
\put(6391,-1366){\makebox(0,0)[lb]{\smash{\SetFigFont{6}{7.2}{\rmdefault}{\mddefault}{\updefault}\special{ps: gsave 0 0 0 setrgbcolor}$E^{\nu}$\special{ps: grestore}}}}
\put(3601,-1411){\makebox(0,0)[lb]{\smash{\SetFigFont{6}{7.2}{\rmdefault}{\mddefault}{\updefault}\special{ps: gsave 0 0 0 setrgbcolor}$E$\special{ps: grestore}}}}
\end{picture}

%% file: exempleincidence.pstex_t
\begin{picture}(0,0)%
\includegraphics{exempleincidence.pstex}%
\end{picture}%
\setlength{\unitlength}{2072sp}%
\begingroup\makeatletter\ifx\SetFigFont\undefined%
\gdef\SetFigFont#1#2#3#4#5{%
  \reset@font\fontsize{#1}{#2pt}%
  \fontfamily{#3}\fontseries{#4}\fontshape{#5}%
  \selectfont}%
\fi\endgroup%
\begin{picture}(6594,2352)(394,-1726)
\put(2701,-1141){\makebox(0,0)[lb]{\smash{\SetFigFont{6}{7.2}{\rmdefault}{\mddefault}{\updefault}\special{ps: gsave 0 0 0 setrgbcolor}The system $S$ on $E$.\special{ps: grestore}}}}
\put(6121,-1366){\makebox(0,0)[lb]{\smash{\SetFigFont{6}{7.2}{\rmdefault}{\mddefault}{\updefault}\special{ps: gsave 0 0 0 setrgbcolor}$S(E)$\special{ps: grestore}}}}
\put(451,-1366){\makebox(0,0)[lb]{\smash{\SetFigFont{6}{7.2}{\rmdefault}{\mddefault}{\updefault}\special{ps: gsave 0 0 0 setrgbcolor}$E$\special{ps: grestore}}}}
\put(2701,-1411){\makebox(0,0)[lb]{\smash{\SetFigFont{6}{7.2}{\rmdefault}{\mddefault}{\updefault}\special{ps: gsave 0 0 0 setrgbcolor}(Arrows of the form $(p,p)$ .\special{ps: grestore}}}}
\put(2746,-1726){\makebox(0,0)[lb]{\smash{\SetFigFont{6}{7.2}{\rmdefault}{\mddefault}{\updefault}\special{ps: gsave 0 0 0 setrgbcolor}have not been drawn)\special{ps: grestore}}}}
\end{picture}

%% file: stratesG12321l.pstex_t
\begin{picture}(0,0)%
\includegraphics{stratesG12321l.pstex}%
\end{picture}%
\setlength{\unitlength}{2072sp}%
\begingroup\makeatletter\ifx\SetFigFont\undefined%
\gdef\SetFigFont#1#2#3#4#5{%
  \reset@font\fontsize{#1}{#2pt}%
  \fontfamily{#3}\fontseries{#4}\fontshape{#5}%
  \selectfont}%
\fi\endgroup%
\begin{picture}(8079,4929)(484,-4303)
\put(8326,-736){\makebox(0,0)[lb]{\smash{\SetFigFont{6}{7.2}{\rmdefault}{\mddefault}{\updefault}\special{ps: gsave 0 0 0 setrgbcolor}$h$\special{ps: grestore}}}}
\put(901,-1636){\makebox(0,0)[lb]{\smash{\SetFigFont{6}{7.2}{\rmdefault}{\mddefault}{\updefault}\special{ps: gsave 0 0 0 setrgbcolor}$E_{gen}$\special{ps: grestore}}}}
\put(2251,-511){\makebox(0,0)[lb]{\smash{\SetFigFont{6}{7.2}{\rmdefault}{\mddefault}{\updefault}\special{ps: gsave 0 0 0 setrgbcolor}$a$\special{ps: grestore}}}}
\put(2251,-2986){\makebox(0,0)[lb]{\smash{\SetFigFont{6}{7.2}{\rmdefault}{\mddefault}{\updefault}\special{ps: gsave 0 0 0 setrgbcolor}$b$\special{ps: grestore}}}}
\put(4726,164){\makebox(0,0)[lb]{\smash{\SetFigFont{6}{7.2}{\rmdefault}{\mddefault}{\updefault}\special{ps: gsave 0 0 0 setrgbcolor}$c$\special{ps: grestore}}}}
\put(3826,-1186){\makebox(0,0)[lb]{\smash{\SetFigFont{6}{7.2}{\rmdefault}{\mddefault}{\updefault}\special{ps: gsave 0 0 0 setrgbcolor}$d$\special{ps: grestore}}}}
\put(4501,-2986){\makebox(0,0)[lb]{\smash{\SetFigFont{6}{7.2}{\rmdefault}{\mddefault}{\updefault}\special{ps: gsave 0 0 0 setrgbcolor}$e$\special{ps: grestore}}}}
\put(6976,-736){\makebox(0,0)[lb]{\smash{\SetFigFont{6}{7.2}{\rmdefault}{\mddefault}{\updefault}\special{ps: gsave 0 0 0 setrgbcolor}$f$\special{ps: grestore}}}}
\put(6976,-1861){\makebox(0,0)[lb]{\smash{\SetFigFont{6}{7.2}{\rmdefault}{\mddefault}{\updefault}\special{ps: gsave 0 0 0 setrgbcolor}$g$\special{ps: grestore}}}}
\end{picture}

%% file: systflechesl.pstex_t
\begin{picture}(0,0)%
\includegraphics{systflechesl.pstex}%
\end{picture}%
\setlength{\unitlength}{1579sp}%
\begingroup\makeatletter\ifx\SetFigFont\undefined%
\gdef\SetFigFont#1#2#3#4#5{%
  \reset@font\fontsize{#1}{#2pt}%
  \fontfamily{#3}\fontseries{#4}\fontshape{#5}%
  \selectfont}%
\fi\endgroup%
\begin{picture}(13225,2976)(2239,-2592)
\put(3911,-2516){\makebox(0,0)[lb]{\smash{\SetFigFont{7}{8.4}{\rmdefault}{\mddefault}{\updefault}\special{ps: gsave 0 0 0 setrgbcolor}Arrow Systems(A point symbolizes an arrow from a monomial to itself)\special{ps: grestore}}}}
\put(5732,-1966){\makebox(0,0)[lb]{\smash{\SetFigFont{7}{8.4}{\rmdefault}{\mddefault}{\updefault}\special{ps: gsave 0 0 0 setrgbcolor}From $b$ to $d$\special{ps: grestore}}}}
\put(11297,-1936){\makebox(0,0)[lb]{\smash{\SetFigFont{7}{8.4}{\rmdefault}{\mddefault}{\updefault}\special{ps: gsave 0 0 0 setrgbcolor}From $c$ to $e$\special{ps: grestore}}}}
\put(13847,-1951){\makebox(0,0)[lb]{\smash{\SetFigFont{7}{8.4}{\rmdefault}{\mddefault}{\updefault}\special{ps: gsave 0 0 0 setrgbcolor}From $d$ to $f$\special{ps: grestore}}}}
\put(8552,-1951){\makebox(0,0)[lb]{\smash{\SetFigFont{7}{8.4}{\rmdefault}{\mddefault}{\updefault}\special{ps: gsave 0 0 0 setrgbcolor}From $b$ to $f$\special{ps: grestore}}}}
\put(2701,-2011){\makebox(0,0)[lb]{\smash{\SetFigFont{7}{8.4}{\rmdefault}{\mddefault}{\updefault}\special{ps: gsave 0 0 0 setrgbcolor}From $a$ to $e$\special{ps: grestore}}}}
\end{picture}